\documentclass{sig-alternate}
\usepackage{comment}
\usepackage{tikz}
\usepackage{xspace}
\usepackage{url}
\usepackage{xspace}

\newcommand{\bit}{\begin{itemize}}
\newcommand{\eit}{\end{itemize}}
\newcommand{\ben}{\begin{enumerate}}
\newcommand{\een}{\end{enumerate}}
\newcommand{\beq}{\begin{equation}}
\newcommand{\eeq}{\end{equation}}

\newcommand{\reals}{{\mbox{\bf R}}}

\newcommand{\cvxjl}{\texttt{Convex}\xspace} 

\newcommand{\lambdamin}{\lambda_{\rm min}}

\newcommand{\etal}{{\it et. al.}\xspace} 
\newcommand{\ie}{{\it i.e.}}

\begin{document}
\conferenceinfo{HPTCDL}{November 16-21, 2014, New Orleans, Louisiana, USA}
\CopyrightYear{2014}
\crdata{978-1-4799-5500-8/14}

\title{Convex Optimization in Julia}

\numberofauthors{6}
\author{
\alignauthor Madeleine Udell\\
       \email{udell@stanford.edu}
\alignauthor Karanveer Mohan\\
       \email{kvmohan@stanford.edu}
\alignauthor David Zeng\\
       \email{dzeng0@stanford.edu}
\and  
\alignauthor Jenny Hong\\
       \email{jyunhong@stanford.edu}
\alignauthor Steven Diamond\\
       \email{stevend2@stanford.edu}
\alignauthor Stephen Boyd\\
       \email{boyd@stanford.edu}
}

\maketitle
\begin{abstract}
This paper describes \cvxjl\footnote{available online at \url{http://github.com/cvxgrp/Convex.jl}}, 
a convex optimization modeling framework in Julia.
\cvxjl translates problems from a user-friendly functional language into an abstract syntax tree 
describing the problem. This concise representation of the global structure of the problem 
allows \cvxjl to infer whether the problem complies with the rules of disciplined convex programming (DCP),
and to pass the problem to a suitable solver. 
These operations are carried out in Julia using multiple dispatch, 
which dramatically reduces the time
required to verify DCP compliance and to parse a problem into conic form. 
\cvxjl then automatically chooses an appropriate backend solver to solve the conic form problem.
\end{abstract}

\category{G.1.6}{Numerical Analysis}{Optimization}[Convex Programming]

\category{D.2.4}{Software Engineering}{Software/Program
Verification}[Model checking]

\terms{Languages}

\keywords{Convex programming, automatic verification, symbolic computation,
multiple dispatch}

\section{Introduction}

The purpose of an optimization modeling language is to translate an
optimization problem from a user-friendly language into a solver-friendly
language. In this paper, we present an approach to bridging the gap. We show
that mathematical optimization problems can be parsed from a simple
human-readable form, using the trope of function composition, into an
abstract syntax tree (AST) representing the problem.

Representing the problem as an AST facilitates subsequent computations.
For example, \cvxjl can efficiently check if a problem is convex by applying the rules of
\emph{disciplined convex programming} (DCP), pioneered by Michael Grant and Stephen
Boyd in \cite{CVX, GBY:06}.  \cvxjl can also use the AST to convert the problem into a
conic form optimization problem, allowing a solver access to a complete and
computationally concise global description of the problem \cite{nesterov2004}.

\paragraph{Julia}

The Julia language \cite{bezanson2012} is a high-level, high-performance
dynamic programming language for technical computing. With a syntax familiar to
users of other technical computing languages such as Matlab, it takes advantage
of LLVM-based just-in-time (JIT) compilation \cite{llvm} to approach and often
match the performance of C \cite{julialang}.  In \cvxjl, we make particular use
of \emph{multiple dispatch} in Julia, an object-oriented paradigm in which
different methods may be called to implement a function
depending on the data types (classes) of the arguments to the function, 
rather than having functions encapsulated inside classes \cite{bezanson2014}.

The \cvxjl project supports the assertion that ``technical computing is [the]
killer application [of multiple dispatch]'' \cite{bezanson2014}.  We show in
this paper that Julia's multiple dispatch allows the authors of a technical
computing package to write extremely performant code using a high level of
abstraction.  Indeed, the abstraction and generality of the code has pushed the
authors of this paper toward more abstract and general formulations of the
mathematics of disciplined convex programming,
while producing code whose performance rivals and often surpasses codes with
similar functionality in other languages.

Moreover, multiple dispatch enforces a separation of the fundamental objects in
\cvxjl --- functions --- from the methods for operating on them. This
separation has other benefits: in \cvxjl, it is easy to implement new methods
for operating on the AST of an optimization problem, opening the door to new
structural paradigms for parsing and solving optimization problems.

\paragraph{Mathematical optimization}

A traditional and familiar form for specifying a mathematical optimization
problem (MOP) is to write
\[
\begin{array}{ll}
\mbox{minimize} & f_0(x) \\
\mbox{subject to} & f_i(x) \leq 0, \quad i=1,\ldots,m \\
& h_i(x) = 0, \quad i=1,\ldots,p,
\end{array}
\]
with variable $x\in \reals^n$.
Here, a problem instance is specified as list of \emph{functions}
$f_0, \ldots, f_m, h_1, \ldots, h_p: \reals^n \to \reals$.
The function $f_0$ is called the objective function;
$f_1 \ldots, f_m$ are called the inequality constraint functions; and
$h_1, \ldots, h_p$ are called the equality constraint functions.
The field of mathematical optimization is concerned with finding methods for
solving MOP.

Some special structural forms of MOP are well known. For example, if all the
functions $f_0, \ldots, f_m$ and $h_1, \ldots, h_p$ are affine, the MOP is called a
linear program (LP).  When each of the functions $f_0, \ldots, f_m$ is convex,
and $h_1, \ldots, h_p$ are all affine, it is called a convex program (CP).
This paper focusses on solving CPs, which can be solved much more quickly than
general MOPs \cite{johnson2010}.

Perhaps the most consistent theme of the convex optimization literature
concerns the importance of using the \emph{structural form} of the CP (that is,
the properties of the functions $f_0, \ldots, f_m, h_1, \ldots, h_p$) in order
to devise faster solution methods.  The advantage of this structural approach
is that it allows for a division of labor between \emph{users} of optimization
and \emph{designers} of optimization algorithms.  Designers of optimization
algorithms write solvers specialized to a particular structural form of MOP.
Users can take advantage of innovations in optimization algorithms and theory
so long as they can \emph{identify} the structural form of their problem, and
apply the appropriate solver.

A number of general purpose solvers have been developed, including ACADO
\cite{HFD:11,ACADO,ACADOMatlab}, Ipopt \cite{ipopt}, NLopt
\cite{johnson2010}), and the solvers in AMPL
\cite{fourer1993}, GAMS \cite{brooke1998,bussieck2004}, and Excel \cite{ExcelSolver},
along with associated modeling languages that make it easy
to formulate problems for these solvers, including AMPL \cite{fourer1993}, GAMS
\cite{brooke1998,bussieck2004}, and JuMP \cite{lubin2013}.  However, these
solvers can be significantly slower than ones designed for special problem
classes \cite{karmarkar1984,nemirovsky1983,nesterov1997}. A modeling language
that retains information about the structure of the problem can often choose a solver
that will significantly outperform a general purpose solver.

\paragraph{Cone programs}
\cvxjl targets a special class of convex programs called \emph{cone programs}, which
includes LPs, QPs, QCQPs, SOCPs, SDPs, and exponential cone programs as special
cases (see \S\ref{s-cone} for more details).  Cone programs occur frequently in
a variety of fields \cite{LVBL:98, vandenberghe1996, wolkowicz2000}, can be
solved quickly, both in theory and in practice \cite{NN:94}.

A number of modeling systems have been developed to automatically perform
translations of a problem into standard conic form.  The first parser-solver,
SDPSOL, was written using Bison and Flex in 1996, and was able to automatically
parse and solve SDPs with a few hundred variables \cite{WB:95}.
YALMIP \cite{Lofberg:04} was its first modern successor, followed shortly
thereafter by CVX \cite{CVX}, both of which embedded a convex optimization
modeling language into MATLAB, a proprietary language.  Diamond \etal followed
with CVXPY~\cite{cvxpy}, a convex optimization modeling language in python
which uses an object-oriented approach.  The present paper concerns \cvxjl,
which borrows many ideas from CVXPY, but takes advantage of language features
in Julia (notably, multiple dispatch) to make the modeling layer simpler and
faster.

A few more specialized conic parsing/modeling tools have also appeared
targeting small QPs and SOCPs with very tight performance requirements,
including CVXGEN \cite{JB:12} and QCML \cite{chu2013}.

\paragraph{Organization}
This paper is organized as follows.  In
\S\ref{s-representation}, we discuss how \cvxjl represents optimization
problems, and show how this representation suits the requirements of both
users and solvers.  In \S\ref{s-dcp}, we show how to use the AST to verify
convexity of the problem using the rules of disciplined convex programming.  In
\S\ref{s-cone}, we show how to use the AST to transform the problem into conic
form, and pass the problem to a conic solver.  Finally, in \S\ref{s-speed} we
compare the performance of \cvxjl to other codes for disciplined convex
programming.

\section{Functional representation}\label{s-representation}
In this section, we describe the functional representation used by \cvxjl.
The basic building block of \cvxjl is called an \emph{expression},
which can represent a variable, a constant, or a function of another expression.
We discuss each kind of expression in turn.

\paragraph{Variables}

The simplest kind of expression in \cvxjl is a variable.
Variables in \cvxjl are declared using the \verb|Variable|
keyword, along with the dimensions of the variable.

\begin{verbatim}
# Scalar variable
x = Variable()

# 4x1 vector variable
y = Variable(4)

# 4x2 matrix variable
z = Variable(4, 2)
\end{verbatim}

Variables may also be declared as having special properties, such as being
(entrywise) positive or negative, or (for a matrix) being symmetric, with 
nonnegative eigenvalues (\ie, positive semidefinite).
These properties are treated as constraints in any problem constructed using
the variables.

\begin{verbatim}
# Positive scalar variable
x = Variable(Positive())

# Negative 4x1 vector variable
y = Variable(4, Negative())

# Symmetric positive semidefinite
# 4x4 matrix variable
z = Semidefinite(4)
\end{verbatim}


\paragraph{Constants}

Numbers, vectors, and matrices present in the Julia environment are wrapped
automatically into a \verb|Constant| expression when used in a \cvxjl expression.  
The \verb|Constant| expression stores a pointer to the underlying value, so the
expression and any derived expression will reflect changes made to the value 
after the expression is formed.  
In other words, a \cvxjl expression is parametrized by the constants it encapsulates.

\paragraph{Atoms}

Functions of expressions may be composed using elements called \emph{atoms},
which represent simple functions with known properties such as
\verb|+|, \verb|*|, \verb|abs|, and \verb|norm|.  
A list of some of the atoms available in \cvxjl is given in Table
\ref{t-atoms}.

Atoms are applied to expressions using operator overloading. Hence, \verb|2+2|
calls Julia's built-in addition operator, while \verb|2+x| calls the \cvxjl
addition method and returns a \cvxjl expression with top-level atom \verb|+| 
and arguments \verb|2| and \verb|x|.  
Many of the useful language
features in Julia, such as arithmetic, array indexing, and matrix transpose are
overloaded in \cvxjl so they may be used with variables and expressions 
just as they are used with native Julia types.

For example, all of the following form valid expressions.
\begin{verbatim}
# indexing, multiplication, addition
e1 = y[1] + 2*x

# expressions can be affine, convex, or concave
e3 = sqrt(x) + log(x)

# more atoms
e2 = 4 * pos(x) + max(abs(y)) + norm(z[:,1],2)
\end{verbatim}

For each atom, we define five methods that allow \cvxjl to solve problems involving
that atom. The methods are called on expressions, and are dispatched 
based on the type (\ie, the top-level atom) of the expression.
\bit
\item
\verb|sign|
returns the \emph{sign} of the function the top-level atom represents
(positive, negative, or no sign) over the range of its arguments.
\item
\verb|evaluate| 
evaluates the function the top-level atom represents on the values of its arguments.
\item
\verb|curvature| 
returns the curvature of the function the top-level atom represents (convex, concave,
affine, or neither) over the range of its arguments.
\item
\verb|monotonicity| 
returns the curvature of the function the atom represents 
(nondecreasing, nonincreasing, or not monotonic) 
over the range of its arguments.
\item
\verb|conic_form| returns the conic form of an expression; see \S\ref{s-cone} for details.
\eit

These methods are computed 
using local combination rules,
so they are never wrong, 
but sometimes miss strong classifications that cannot be deduced from local information.  
For example, \verb|square(x) + x - x| is nonnegative for any value of x, but \cvxjl will return no sign.

To define a new atom that can be used to solve problems in \cvxjl,
it suffices to define these five methods.
We discuss how these methods are used to verify convexity and solve conic problems 
in \S\ref{s-dcp}) and \S\ref{s-cone}, respectively.

Atoms may also be defined by composing other atoms; we call these \emph{composite}
atoms.  For example, \verb|norm(x, 1)| is implemented as \verb|sum(abs(x))|.

\paragraph{Expressions}

Formally, expressions are defined recursively to be a variable, a constant, or
the result of applying an atom to one or more expressions.  In \cvxjl, each
expression is represented as an AST with a variable or constant at each leaf
node, and a atom at each internal node. 
We refer to the arguments to the atom at a certain node as the \emph{children} of that node.  
The top-level atom is called the \emph{head} of the AST.

Considered as a tree without leaves, the AST may be seen as representing
a function formed by composing the atoms at each node according to the structure of the tree; 
with the leaves, it may be interpreted as a function with a closure, a curried
function, or simply a function of named arguments.  

The AST is a directed acyclic graph (DAG) \cite{FM:10}. 
The idea of using a DAG to study global
properties of optimization problems dates back to work by Kantorovich in the
1950s \cite{kantorovich1957}.  

Expressions may be evaluated when every leaf node in their AST has
already been assigned a value.  For example,
\begin{verbatim}
x = Variable()
e = max(x,0)
x.value = -4
evaluate!(e)
\end{verbatim}
evaluates to \verb|0|.

One novel aspect of \cvxjl is that expressions are identified via a unique identifier
computed as a hash on the name of the atom at the head of the AST, together with
the unique identifiers of its children in the AST.  
Thus, two different instances of the same expression will
automatically be assigned the same unique identifier.  This identification allows \cvxjl
to identify when the same expression is used in multiple places in the problem,
to memoize and reuse the results of computations on a given expression,
and to reduce the complexity of the problem that is ultimately passed to the solver.

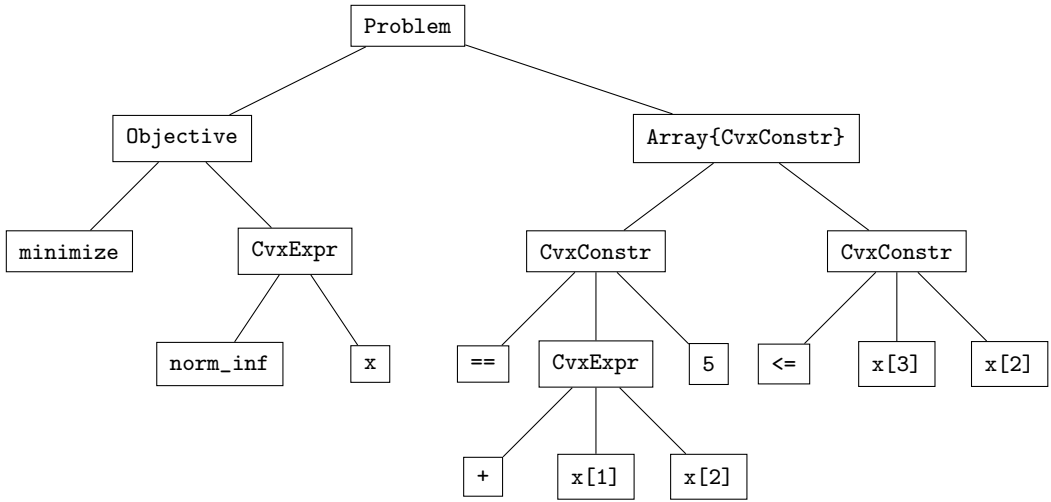
\begin{figure*}
\begin{center}
\begin{tikzpicture}
[dot/.style={rectangle,draw=black,fill=white,inner sep=5pt,minimum size=5pt}]
\node[dot] at (-2,10) (n1) {\texttt{Problem}};
\node[dot] at (-5,8.5) (n2) {\texttt{Objective}};
\node[dot] at (2.5,8.5) (n3) {\texttt{Array\{CvxConstr\}}};
\node[dot] at (-6.5,7) (n4) {\texttt{minimize}};
\node[dot] at (-3.5,7) (n5) {\texttt{CvxExpr}};
\node[dot] at (0.5,7) (n6) {\texttt{CvxConstr}};
\node[dot] at (4.5,7) (n7) {\texttt{CvxConstr}};
\node[dot] at (-4.5,5.5) (n8) {\texttt{norm\_inf}};
\node[dot] at (-2.5,5.5) (n9) {\texttt{x}};
\node[dot] at (-1,5.5) (n10) {\texttt{==}};
\node[dot] at (0.5,5.5) (n11) {\texttt{CvxExpr}};
\node[dot] at (2,5.5) (n12) {\texttt{5}};
\node[dot] at (3,5.5) (n13) {\texttt{<=}};
\node[dot] at (4.5,5.5) (n14) {\texttt{x[3]}};
\node[dot] at (6,5.5) (n15) {\texttt{x[2]}};
\node[dot] at (-1,4) (n16) {\texttt{+}};
\node[dot] at (0.5,4) (n17) {\texttt{x[1]}};
\node[dot] at (2,4) (n18) {\texttt{x[2]}};

\draw[-] (n1) -- (n2);
\draw[-] (n1) -- (n3);
\draw[-] (n2) -- (n4);
\draw[-] (n2) -- (n5);
\draw[-] (n3) -- (n6);
\draw[-] (n3) -- (n7);
\draw[-] (n5) -- (n8);
\draw[-] (n5) -- (n9);
\draw[-] (n6) -- (n10);
\draw[-] (n6) -- (n11);
\draw[-] (n6) -- (n12);
\draw[-] (n7) -- (n13);
\draw[-] (n7) -- (n14);
\draw[-] (n7) -- (n15);
\draw[-] (n11) -- (n16);
\draw[-] (n11) -- (n17);
\draw[-] (n11) -- (n18);
\end{tikzpicture}
\end{center}
\caption{Graphical representation of a problem}
\label{f-expr-tree}
\end{figure*}

\paragraph{Constraints}

\emph{Constraints} in \cvxjl are declared using the standard comparison
operators \verb|<=|, \verb|>=|, and \verb|==|.  They specify relations that
must hold between two expressions.  \cvxjl does not distinguish between strict
and non-strict inequality constraints.

\begin{verbatim}
# affine equality constraint
A = randn(3,4); b = randn(3)
c1 = A*y == b

# convex inequality constraint
c2 = norm(y,2) <= 2
\end{verbatim}


\paragraph{Problems}
A \emph{problem} in \cvxjl consists of a \emph{sense} (minimize, maximize,
or satisfy), an objective (an expression to which the sense verb is to be
applied), and a list of zero or more constraints which must be satisfied at the
solution.

For example, the problem
\[
\begin{array}{ll}
\mbox{minimize}   & \| x \|_\infty \\
\mbox{subject to} & x_1 + x_2 = 5 \\
& x_3 \leq x_2,
\end{array}
\]
where $x \in \reals^3$ is the optimization variable, can be expressed in \cvxjl as
\begin{verbatim}
x = Variable(3)
constraints = [x[1]+x[2] == 5, x[3] <= x[2]]
p = minimize(norm_inf(x), constraints)
\end{verbatim}

When the user types this code, the input problem description is parsed
internally into a \verb|Problem| object with two attributes: an objective and
a list of constraints. Each expression appearing in the problem is
represented as an AST.
For example, the problem \verb|p| constructed above has the following properties:
\begin{verbatim}
p.objective   = (:norm_inf, x)
p.constraints = [(:==, 5, (:+, (:getindex, x, 1),
                             (:getindex, x, 2))),
            (:<=, getindex(x, 3), getindex(x, 2)]
\end{verbatim}
Here, we use Polish (prefix) notation to display the ASTs.
Figure \ref{f-expr-tree} gives a graphical representation of the structure 
of problem \verb|p|.

\paragraph{Solving problems}


The \verb|solve!| method in \cvxjl
checks that the problem is a disciplined convex program, converts it into conic
form, and passes the problem to a solver:
\begin{verbatim}
solve!(p)
\end{verbatim}
After \verb|solve!| has been called on problem \verb|p|, the optimal
value of \verb|p| can be queried with \verb|p.optval|, and any expression
\verb|x| used in \verb|p| is annotated with a value, accessible via \verb|x.value|.
If the solver computes dual variables, these are populated in the constraints;
for example, the dual variable corresponding to the first constraint can be
accessed as \verb|p.constraints[1].dual_value|.

\section{Disciplined convex programming}\label{s-dcp}

Checking if a function is convex is a difficult problem in general.  Many
approaches to verifying convexity have been proposed
\cite{Stout:78,Crusius:03}, and implemented in modeling systems ranging from
AMPL to Microsoft Excel \cite{FM:10,NFK:04,ExcelSolver}.  

In \cvxjl,
we have chosen to use the framework of \emph{disciplined convex programming}
(DCP) to verify problem convexity, which has a number of advantages.  First,
the simplicity of the DCP rules makes the determination of convexity
transparent to the user \cite{BV:04}.
Second, a problem that complies with the DCP rules can be converted
to an equivalent conic form problem.  Hence \cvxjl is able to take advantage of
extremely fast and reliable solvers for any problem whose convexity it can verify.

\paragraph{Disciplined convex expressions}
In \cvxjl, we define five different kinds of curvature: constant, affine, convex, concave and not DCP.
The simplest expressions are single variables, which are affine,
and constants, which are (unsurprisingly) constant.
A constant expression is also affine, and an expression is affine if and only if it is convex and concave.

We determine the curvature of other expressions recursively using the following rule.
Suppose that the curvatures of the expressions $e_1,\ldots,e_n$ are known.
Then the expression $f(e_1,\ldots,e_n)$ is convex if the function $f$ is convex
on the range of $(e_1,\ldots,e_n)$, and for each $i=1,\ldots,n$,
\bit
\item $f$ is nondecreasing in its $i$th argument over the range of $(e_1,\ldots,e_n)$, and $e_i$ is convex;
\item $f$ is nonincreasing in its $i$th argument over the range of $(e_1,\ldots,e_n)$, and $e_i$ is concave; or
\item $e_i$ is affine.
\eit
The expression $f(e_1,\ldots,e_n)$ is concave if $(-f)(e_1,\ldots,e_n)$ is convex,
and it is affine if it is both convex and concave.
Note that the curvature of an expression need not be the same as the curvature of the 
top-level atom in that expression.
For example, the atom \verb|+| has affine curvature, but the expression \verb|2 + square(x)| is convex.

Disciplined convex expressions (also called DCP compliant, or simply DCP, expressions)
are defined inductively.
We say that an expression is DCP if its curvature (affine, convex, or
concave) can be inferred from these rules;
otherwise, we say it is not DCP.  
Hence in order to decide if an expression
is DCP, we must be able to calculate the curvature of $f$
and the monotonicity of $f$ in each argument for every atom $f$ appearing in the expression.

The set of DCP expressions depends on the set of atoms. 
For example, $\log(\sum_{i=1}^n \exp{x_i})$ is convex, but
its convexity cannot be derived from the DCP rules using only the atoms
$\log$ and $\exp$ \cite{GBY:06}.  
In \cvxjl, it is easy to add new atoms, 
allowing sophisticated users to expand the set of
functions that \cvxjl can recognize as convex.

The set of DCP expressions also depends on what is known
about the monotonicity of atoms as a function of their arguments.  For
example, the quadratic-over-linear function $f(x,y) = x^T x / y$ is not convex
on $\reals^2$, but it is convex on $\reals \times (0, \infty)$.  This
observation is known as a \emph{signed} convexity rule \cite{cvxpy}. Our
formulation of the DCP rule above in terms of the \emph{ranges} of the
expressions $(e_1, \ldots, e_n)$ generalizes this observation.
\cvxjl implements a signed monotonicity rule: the monotonicity of an atom depends
on the sign of the atom's children.


\paragraph{Verifying DCP using multiple dispatch}

\cvxjl checks that expressions are DCP using multiple dispatch.

Arithmetic is defined directly on \verb|Vexity| and \verb|Monotonicity| types so that addition
of curvatures and multiplication of curvatures by monotonicities enforce the DCP
rule. 
For example, adding \verb|ConvexVexity| to \verb|AffineVexity| results in \verb|ConvexVexity|;
multiplying \verb|ConvexVexity| by \verb|Nonincreasing| monotonicity results in \verb|ConcaveVexity|.
The curvature of an expression can then be calculated by applying the following function:
\begin{verbatim}
function vexity(x::AbstractExpr)
  monotonicities = monotonicity(x)
  vexity = curvature(x)
  for i = 1:length(x.children)
    vexity += monotonicities[i]
                * vexity(x.children[i])
  end
  return vexity
end
\end{verbatim}

This approach to checking convexity has a number of advantages over using
if-else statements to enforce the DCP rules.  First, the implementation is
aesthetically appealing: the code follows the mathematics.  Second, this 
structure makes it easy to write new atoms, 
since one only needs to implement the \verb|curvature| and
\verb|monotonicity| methods to allow \cvxjl to determine the curvature of expressions
involving the new atom.  Third, it reduces the time needed to verify
convexity, since multiple dispatch is implemented as a lookup table rather than
via a (slower) \verb|if| statement.

\paragraph{Disciplined convex constraints}

A constraint is called a \emph{disciplined convex constraint} (or DCP constraint) if it has the form
\bit
\item $e_1 \leq e_2$, where $e_1$ is a convex DCP expression and $e_2$ is a concave DCP expression;
\item $e_1 \geq e_2$, where $e_1$ is a concave DCP expression and $e_2$ is a convex DCP expression; or
\item $e_1 = e_2$, where $e_1$ and $e_2$ are both affine DCP expressions.
\eit
(For the purposes of this definition, remember that constant expressions are affine, and affine
expressions are both convex and concave.)

\paragraph{Disciplined convex programs}
A problem is called a \emph{disciplined convex program} if
\ben
\item the (sense, objective) are
\bit
\item (minimize, convex DCP expression),
\item (maximize, concave DCP expression), or
\item (satisfy,); and
\eit
\item every constraint in the problem is DCP compliant.
\een
(Note that the sense satisfy does not take an objective.) 
When an optimization problem is formed, \cvxjl applies the DCP rules recursively to 
determine whether the problem is DCP.

\section{Conic form optimization}\label{s-cone}

A conic form optimization problem is written as
\beq
\begin{array}{ll}
\mbox{minimize} & c^T x \\
\mbox{subject to} & A x = b \\
& x \in \mathcal K,
\end{array}
\label{eq-conic-problem}
\eeq
with variable $x$.
Here $\mathcal K$ is a \emph{cone}: a set of points such that
$x \in \mathcal K$ iff $rx \in \mathcal K$ for every $r\geq 0$.
If in addition $\mathcal K$ is convex, it is called a convex cone.

Examples of convex cones include the following.
\begin{itemize}
\item \emph{The zero cone.} $\mathcal K^1_0 = \{0\}$
\item \emph{The free cone.} $\mathcal K^n_{\mathrm{free}} = \reals^n$
\item \emph{The positive orthant.} $\mathcal K^n_+ = \{x \in \reals^n: x\geq 0\}$
\item \emph{The second order cone.} \\
$\mathcal K^{n+1}_{\mathrm{SOC}} = \{(x,t) \in \reals^n: \|x\| \leq t\}$
\item \emph{The semidefinite cone.} \\
$\mathcal K^{n^2}_{\mathrm{SDP}} = \{X \in \reals^{n\times n}: \lambdamin(X) \geq 0, X=X^T\}$
\item \emph{The exponential cone.} \\
$\mathcal K^{3}_{\mathrm{exp}} = \{(x,y,z) \in \reals^3 : y e^{x/y} \leq z, y>0 \}$
\end{itemize}
In general, a cone $\mathcal K$ may also be given as a product of simpler cones,
\[
\mathcal K = \mathcal K_1 \times \cdots \times \mathcal K_p.
\]

A conic form optimization problem is specified by the problem data $A$, $b$,
$c$, and a cone $K$. 
In this paper, we extend this definition: we say a problem is in conic form 
if every expression in the problem is affine, 
and the constraints are all conic constraints,
since it is trivial to rewrite such a problem in the form \ref{eq-conic-problem}.

A number of fast solvers exist for conic form problems in a variety of languages. 
We list some of these solvers in Table~\ref{t-solvers}.

\cvxjl is able to automatically rewrite a problem specified as a DCP in
standard cone format, and to pick a cone solver automatically depending on the
cones required to represent the problem and the solvers available on the user's computer,
using the \verb|loadconicproblem!| interface in MathProgBase \cite{lubin2013}.
This allows \cvxjl to, for example, choose a SOCP solver if no exponential or
SDP cones are present in the problem, which can result in significantly faster
solve time \cite{alizadeh2003}.

\subsection{Example}
Here, we describe how \cvxjl transforms a DCP problem into conic form.

\paragraph{Conic form expression}
The value of the function $f(x)$ (for fixed $x$) is the optimal value of the
trivial optimization problem 
\[
\mbox{minimize} \quad f(x),
\]
with no variable. 
We say a function is \emph{cone-representable} if the value of the function
is the optimal value of a (nontrival) conic form optimization problem.
This conic form optimization problem is usually constructed by introducing auxiliary variables.
For example, $f(x) = |x|$ is cone-representable, since $|x|$ is the optimal value of the 
conic problem
\[
\begin{array}{ll}
\mbox{minimize} & t \\
\mbox{subject to} & t - x \in \mathcal K_+ \\
& t + x \in \mathcal K_+,
\end{array}
\]
with variable $t$.

In general, define the conic form template of a function $f$ to be a constructor for a
conic form optimization problem.  The conic form template of $f$ takes a list of arguments 
$x_1, \ldots,x_n$, and returns a conic form optimization problem whose optimal value
is $f(x_1,\ldots,x_n)$.
Concretely, it returns a sense, an objective, and a (possibly empty) set of
constraints $(s, o, \mathcal C)$ such that every expression appearing in the objective 
or constraints is affine in $x_1,\ldots, x_n$.
Note that the objective of the conic form problem may be vector or matrix
valued.  For example, the conic form template for $f(x) = Ax$ is simply 
\[
\mbox{minimize} \quad A x
\]
(as a function of $x$).

The conic form template of a function is also known as the
\emph{graph form} of the function \cite{GB:08}.  For convex functions $f$, any
feasible point $(x_1, \ldots, x_n,t)$ for the conic form problem is a point in
the epigraph $\{(x,t): f(x) \leq t\}$ of $f$; for concave functions, any
feasible point of the conic form problem lies in the hypograph of $f$.  

\cvxjl uses the conic form templates of the atoms in an expression recursively
in order to construct the conic form of a expression.  Recall that
every expression is represented as an AST with a function (atom) as its top-level node.  We
refer to the arguments to the function at a certain node as the children of
that node.  The conic form of an expression $e$ is constructed
recursively as follows:
\bit
\item If $e$ is affine, return $(e,\emptyset)$.
\item Otherwise:
\ben
\item Compute the conic forms $(s_i, o_i, \mathcal C_i)$ 
for each child $i=1,\ldots,n$ of the top-level node of $e$.
\item Apply the conic form template for the top-level node to the list of objectives
$(o_1,\ldots,o_n)$, 
producing a conic form problem $(s, o, \mathcal C)$.
\item Return $\left( s, o, \mathcal C \cup \left(\cup_i \mathcal C_i\right)\right)$.
\een
\eit


\paragraph{Conic form problem}

Now that we have defined the conic form
for an expression, it is simple to define the conic form
for an optimization problem.
The conic form of an optimization problem, given as an objective
and constraint set, is computed as follows:
\ben
\item Let $(s,o,\mathcal C)$ be the conic form of the objective.
\item For each constraint in the constraint set,
\ben
\item Compute the conic forms $(s_l, o_l, \mathcal C_l)$ and $(s_r, o_r, \mathcal C_r)$
for the left and right hand sides of the constraint.
\item Add $\mathcal C_l$ and $\mathcal C_r$ to $\mathcal C$.
\item If the sense of the constraint is
\bit
\item $\leq$, add $o_r - o_l \in \mathcal K_+$ to $\mathcal C$;
\item $\geq$, add $o_l - o_r \in \mathcal K_+$ to $\mathcal C$;
\item $=$, add $o_r - o_l \in \mathcal K_0$ to $\mathcal C$.
\eit
\item Return $( s, o, \mathcal C)$.
\een
\een

Notice that we have not used the sense of the conic forms of the arguments to
an expression in constructing its conic form. This would be worrisome, were it
not for the following theorem:

\newtheorem{theorem}{Theorem}
\begin{theorem}[\cite{GBY:06}]
Let $p$ be a problem with variable $x$ and dual variable $\lambda$,
and let $\Phi(p)$ be the conic form of the problem
with variables $x$ and $t$ and dual variables $\lambda$ and $\mu$.
Here we suppose $t$ and $\mu$ are the primal and dual variables, respectively, that have been introduced in the 
transformation to conic form.
If $p$ is DCP, then any primal-dual solution $(x,t,\lambda,\mu)$ to $\Phi(p)$
provides a solution $(x,\lambda)$ to $p$.
\end{theorem}

To build intuition for this theorem, note that in a DCP expression, a convex
function $f$ will have objectives with the sense $\mbox{minimize}$ spliced into
argument slots in which $f$ is increasing, and objectives with the sense
$\mbox{maximize}$ spliced into argument slots in which $f$ is decreasing.  At
the solution, the coincidence of these senses, monotonicity, and curvature will
force the variables participating in the conic form of each atom to lie
\emph{on} the graph of the atom, rather than simply in the \emph{epigraph}.
Using this reasoning, the theorem can be proved by induction on the depth of
the AST.


\subsection{Solvers}

Julia has a rich ecosystem of optimization routines. The \verb|JuliaOpt|
project collects mathematical optimization solvers, and interfaces to solvers
written in other languages, in a single GitHub repository, while
\verb|MathProgBase| enforces consistent interfaces to these solvers.  Through
integration with \verb|MathProgBase|, \cvxjl is able to use all of the solvers
that accept problems in conic form, which includes all the linear programming
solvers in \verb|JuliaOpt| (including GLPK \cite{glpk} and the commercial
solvers CPLEX \cite{cplex}, Gurobi \cite{cplex}, and Mosek \cite{mosek}), as
well as the open source interior-point SOCP solver ECOS \cite{DB:12}, and the
open source first-order primal-dual conic solver SCS \cite{scs}.  A list of
solvers which can be used with \cvxjl is presented in Table~\ref{t-solvers}.

\section{Speed}\label{s-speed}

Here we present a comparison of \cvxjl with CVXPY\footnote{
available at \url{https://pypi.python.org/pypi/cvxpy} as of \today} \cite{cvxpy}
 and CVX\footnote{
CVX v3.0 beta, available at \url{www.cvxr.com} as of \today} \cite{CVX} 
on a number of representative problems. The problems are chosen to be easy to \emph{solve}
(indeed, they are all easily solved by inspection) but difficult to \emph{parse}. 
We concentrate on problems involving a large number of affine expressions;
other atoms are treated very similarly in the different frameworks.
The code for the tests in the three languages may be found in Appendix~\ref{a-speed}.

We compare both the time required to convert the problem to conic form, and the time required to solve
the resulting problem, since in general the conic form problems produced in
different modeling frameworks may be different.  We use all three modeling
languages with the solver ECOS \cite{DB:12}, and call the solver with the same
parameters (\verb|abstol=1e-7|, \verb|reltol=1e-7|, \verb|feastol=1e-7|, \verb|maxit=100|) 
so that solve times are comparable between the modeling frameworks.

Parse times are presented in Table \ref{t-parse} and solve times in Table \ref{t-solve}.
Parse times are computed by subtracting solve time from the total time to form and solve the problem.  

\begin{table*}
\caption{\label{t-parse} Speed comparisons: parse time (s)}
\begin{center}
\begin{tabular}{|l|| l | l | l | l |} \hline
Test              & CVX   & CVXPY & \cvxjl & \cvxjl compiled \\ \hline
sum               & 2.29  & 4.45  & 5.46   & 1.94            \\ \hline
index             & 3.62  & 9.69  & 8.40   & 5.78            \\ \hline
transpose         & 1.24  & 0.55  & 3.08   & 0.40            \\ \hline
matrix constraint & 1.38  & 0.54  & 2.49   & 0.34            \\ \hline
\end{tabular}
\end{center}
\end{table*}

\begin{table*}
\caption{\label{t-solve} Speed comparisons: solve time (s)}
\begin{center}
\begin{tabular}{|l|| l | l | l | l |} \hline
Test              & CVX   & CVXPY & \cvxjl & \cvxjl compiled \\ \hline
sum               & 0.01  & 5e-4  & 2e-4   & 1e-4            \\ \hline
index             & 0.01  & 0.08  & 0.076  & 0.061           \\ \hline
transpose         & 0.66  & 1.54  & 14.23  & 11.99           \\ \hline
matrix constraint & 0.66  & 2.20  & 3.44   & 3.29            \\ \hline
\end{tabular}
\end{center}
\end{table*}

We present two different timings for \cvxjl. Julia uses just in time (JIT) compilation;
the code is compiled at the first evaluation, and the faster compiled version is used on subsequent calls,
so the first time a function is called is in general slower than the second time.
The third column shows the runtime of the first (uncompiled) evaluation, and the fourth column the second (compiled) time.
While the perfomance of \cvxjl is comparable to CVX and CVXPY at first evaluation, it substantially
outperforms upon the second evaluation.
For applications in which convex optimization routines are called in an \emph{inner} loop of a larger optimization
algorithm (for example, in sequential convex programming), users will see the fast performance on every iteration
of the loop after the first.

\section{Discussion}





The parsing and modeling methodology applied here need not be restricted
to conic form problems.
Indeed, there are many other interesting and useful problem structures. In
general, any problem structure that has been identified and for which there are
specialized, fast solvers may be a good candidate for inclusion into a 
modeling framework like \cvxjl.
Such problem structures include integer linear programs (which have been a focus of the operations research community for decades, 
yielding fast solvers in practice \cite{nemhauser1988, wolsey1998}); biconvex problems (which may
be heuristically solved by alternating minimization \cite{udell2014glrm}); sums of prox-capable
functions (which can be efficiently solved by splitting methods like ADMM, even
when the problem includes spectral functions of large matrices \cite{boyd2011, becker2011});
difference-of-convex programs (for which sequential convex programming often
produces successful solutions in practice \cite{lipp2014, van2003}); and
sigmoidal programming problems (which often admit fast solutions via a
specialized branch and bound method \cite{udell2013sp,udell2014db}); to name only a few.  The number of
structural forms that might be useful far exceeds the number that a potential
consumer of optimization might care to remember.

The framework developed in \cvxjl makes it easy to write new transformations on the AST to detect and transform these problem types into a standard format.  Automatic detection of some of the problem types mentioned above is the subject of ongoing research.

Significant work also remains to allow \cvxjl to handle extremely large scale problems.
For example, \cvxjl currently performs all computations serially; developments in 
Julia's parallelism ecosystem, and particularly its nascent threading capabilities, will enable
\cvxjl to leverage parallelism.

\section{Conclusions}
This paper shows that a convex optimization modeling library may be efficiently and simply implemented using multiple dispatch. The software package \cvxjl uses language features of Julia (notably, multiple dispatch and JIT compilation) to parse DCP problems into conic form as fast or faster than other (significantly more mature) codes for modeling convex optimization.

\section*{Acknowledgments}
This work was developed with support from
the National Science Foundation Graduate Research Fellowship program (under Grant No. DGE-1147470),
the Gabilan Stanford Graduate Fellowship,
the Gerald J. Lieberman Fellowship,
and the DARPA X-DATA program.

\begin{table*}
\centering
\caption{\label{t-atoms} Examples of atoms implemented in \cvxjl}
\begin{tabular}{|l||l|l|p{4cm}|l|} \hline
& Function & Atom & Monotonicity & Curvature \\ \hline
\multicolumn{5}{|l|}{Slicing and shaping atoms} \\ \hline
& $f(x, i) = x_i$ & \texttt{getindex(x,i)}
    & nondecreasing & Affine \\ \hline
& $f(x, y) = [x~y]$ &  \texttt{hcat(x,y)}
    & nondecreasing  & Affine \\ \hline
& $f(x, y) = [x, y]$ &  \texttt{vertcat(x,y)}
    & nondecreasing &  Affine \\ \hline
& $f(X) = (X_{11}, \ldots, X_{nn})$ &  \texttt{diag(X)}
    & nondecreasing & Affine \\ \hline
& $f(X) = X^T$ &  \texttt{transpose(X)}
    & nondecreasing & Affine  \\ \hline
\multicolumn{5}{|l|}{Positive orthant atoms} \\ \hline
& $f(x) = \sum_i x_i$ &  \texttt{sum(x)}
    & nondecreasing & Affine  \\ \hline
& $f(x) = (|x_1|, \ldots, |x_n|)$ &  \texttt{abs(x)}
    & nondecreasing & Convex \\ \hline
& $f(x) = \mbox{max}_i x_i$ &  \texttt{max(x)}
    & nondecreasing & Convex  \\ \hline
& $f(x) = \mbox{min}_i x_i$ &  \texttt{min(x)}
    & nondecreasing & Concave \\ \hline
& $f(x) = \mbox{max}(x, 0)$ &  \texttt{pos(x)}
    & nondecreasing & Convex  \\ \hline
& $f(x) = \mbox{max}(-x, 0)$ &  \texttt{neg(x)}
    & nonincreasing & Convex  \\ \hline
& $f(x) = \|x\|_1$ &  \texttt{norm\_1(x)}
    & nondecreasing $x \geq 0$ \newline nonincreasing $x \leq 0$ & Convex  \\ \hline
& $f(x) = \|x\|_{\infty}$ &  \texttt{norm\_inf(x)}
    & nondecreasing & Convex  \\ \hline
\multicolumn{5}{|l|}{Second-order cone atoms} \\ \hline
& $f(x) = \|x\|_2$ &  \texttt{norm\_2(x)}
    & nondecreasing $x \geq 0$ \newline nonincreasing $x \leq 0$ & Convex  \\ \hline
& $f(X) = \|X\|_F$ &  \texttt{norm\_fro(X)}
    & nondecreasing $X \geq 0$ \newline nonincreasing $X \leq 0$ & Convex  \\ \hline
& $f(x) = x^2$&  \texttt{square(x)}
    & nondecreasing $x \geq 0$ \newline nonincreasing $x \geq 0$ & Convex  \\ \hline
& $f(x) = \sqrt{x}$&  \texttt{sqrt(x)}
    & nondecreasing $x \geq 0$& Concave \\ \hline
& $f(x, y) = \sqrt{xy}$ &  \texttt{geo\_mean(x,y)}
    & nondecreasing & Concave \\ \hline
& $f(x, y) = x^T x / y$ &  \texttt{quad\_over\_lin(x,y)}
    & nonincreasing in $y$ for $y>0$
    \newline nondecreasing in $x$ for $x \geq 0$
    \newline nonincreasing in $x$ for $x \leq 0$
    & Convex  \\ \hline
& $f(x) = 1 / x$ &  \texttt{inv\_pos(x)}
    & nonincreasing & Convex  \\ \hline
& $f(x) = \sum_i x_i^2$ &  \texttt{sum\_squares}
    & nondecreasing in $x$ for $x \geq 0$
    \newline nonincreasing in $x$ for $x \leq 0$
    & Convex  \\ \hline
\multicolumn{5}{|l|}{Exponential cone atoms} \\ \hline
& $f(x) = \exp(x)$ &  \texttt{exp(x)}
    & nondecreasing & Convex  \\ \hline
& $f(X) = \log(x)$ &  \texttt{log(x)}
    & nondecreasing $x > 0$ & Concave  \\ \hline
& $f(x) = \log(\sum_i \exp(x_i))$ &  \texttt{logsumexp(x)}
    & nondecreasing & Convex  \\ \hline
\multicolumn{5}{|l|}{SDP atoms} \\ \hline
& $f(X) = \|X\|_2$ &  \texttt{operatornorm(x)}
    &  & Convex  \\ \hline
& $f(X) = \|X\|_*$ &  \texttt{nuclearnorm(x)}
    &  & Convex  \\ \hline
\end{tabular}
\end{table*}

\begin{table*}
\caption{\label{t-solvers} Conic Form Problem Solvers}
\begin{center}
\begin{tabular}{|l||c c c c c p{4.5cm} c|} \hline
Name & Language & LP & SOCP & SDP & Exp & Method & Open source \\
\hline
CLP \cite{clp}       & C++    & X &   &   &   & Simplex                        & X \\ \hline
Gurobi \cite{gurobi} &        & X & X &   &   & Simplex, Interior Point, \ldots&   \\ \hline
GLPK \cite{glpk}     & C      & X &   &   &   & Simplex, Interior Point, \ldots& X \\ \hline
MOSEK \cite{mosek}   &        & X & X & X &   & Simplex, Interior Point, \ldots&   \\ \hline
ECOS \cite{DB:12}    & C      & X & X &   &   & Interior Point                 & X \\ \hline
SDPA \cite{SDPA}     & C++    & X & X & X &   & Interior Point                 & X \\ \hline
SCS \cite{scs}       & C      & X & X & X & X & Primal-Dual Operator Splitting & X \\ \hline
\end{tabular}
\end{center}
\end{table*}

\bibliographystyle{abbrv}
\bibliography{cvxjl}
%
%


\appendix
\section{Speed tests: code}\label{a-speed} 
\subsection{Summation}
\paragraph{\cvxjl}
\begin{verbatim}
n = 10000
@time begin
x = Variable();
e = 0;
for i = 1:n
  e = e + x;
end
p = minimize(norm2(e-1), x >= 0);
solve!(p, ECOS.ECOSMathProgModel())
end
\end{verbatim}

\paragraph{CVX}
\begin{verbatim}
n = 10000;
tic()
cvx_begin
variable x
e = 0
for i=1:n
    e = e + x;
end
minimize norm(e-1, 2)
subject to
    x >= 0;
cvx_end
toc()
\end{verbatim}

\paragraph{CVXPY}
\begin{verbatim}
%%timeit
n = 10000;
x = Variable()
e = 0
for i in range(n):
    e = e + x
p = Problem(Minimize(norm(e-1,2)), [x>=0])
p.solve("ECOS", verbose=True)
\end{verbatim}

\subsection{Indexing}
\paragraph{\cvxjl}
\begin{verbatim}
n = 10000
@time begin
x = Variable(n);
e = 0;
for i = 1:n
  e = e + x[i];
end
p = minimize(norm2(e-1), x >= 0);
solve!(p, ECOS.ECOSMathProgModel())
end
\end{verbatim}

\paragraph{CVX}
\begin{verbatim}
n = 10000;
tic()
cvx_begin
variable x(n)
e = 0;
for i=1:n
    e = e + x(i);
end
minimize norm(e - 1, 2)
subject to
    x >= 0;
cvx_end
toc()
\end{verbatim}

\paragraph{CVXPY}
\begin{verbatim}
%%timeit
n = 10000
x = Variable(n)
e = 0
for i in range(n):
    e += x[i];
p = Problem(Minimize(norm(e-1,2)), [x>=0])
p.solve("ECOS", verbose=True)
\end{verbatim}

\subsection{Transpose}
\paragraph{\cvxjl}
\begin{verbatim}
n = 500
A = randn(n, n);
@time begin
X = Variable(n, n);
p = minimize(vecnorm(X' - A), X[1,1] == 1);
solve!(p, ECOS.ECOSMathProgModel())
end
\end{verbatim}

\paragraph{CVX}
\begin{verbatim}
n = 500;
A = randn(n, n);
tic()
cvx_begin
variable X(n, n);
minimize(norm(transpose(X) - A, 'fro'))
subject to
    X(1,1) == 1;
cvx_end
toc()
\end{verbatim}

\paragraph{CVXPY}
\begin{verbatim}
n = 500
A = numpy.random.randn(n,n)
%%timeit
X = Variable(n,n)
p = Problem(Minimize(norm(X.T-A,'fro')), [X[1,1] == 1])
p.solve("ECOS", verbose=True)
\end{verbatim}

\subsection{Matrix constraint}
\paragraph{\cvxjl}
\begin{verbatim}
n = 500;
A = randn(n, n);
B = randn(n, n);
@time begin
X = Variable(n, n);
p = minimize(vecnorm(X - A), X == B);
solve!(p, ECOS.ECOSMathProgModel())
end
\end{verbatim}

\paragraph{CVX}
\begin{verbatim}
n = 500;
A = randn(n, n);
B = randn(n, n);
tic()
cvx_begin
variable X(n, n);
minimize(norm(X - A, 'fro'))
subject to
    X == B;
cvx_end
toc()
\end{verbatim}

\paragraph{CVXPY}
\begin{verbatim}
n = 500
A = numpy.random.randn(n,n)
B = numpy.random.randn(n,n)
%%timeit
X = Variable(n,n)
p = Problem(Minimize(norm(X-A,'fro')), [X == B])
p.solve(verbose=True)
\end{verbatim}

\end{document}